\documentclass[english, 11pt]{amsart}

\usepackage{tikz}
\usepackage{aliascnt}
\usepackage{amsfonts}
\usepackage{mathrsfs}
\usepackage{amsthm}
\usepackage{amsmath,amssymb}
\usepackage[all,poly,knot]{xy}
\usepackage{amsfonts}
\usepackage{color}
\usepackage{hyperref}
\usepackage{comment}  
\usepackage{csquotes}   
\usepackage{amscd}

\usepackage{contour}
\usepackage{ulem}

\newcommand{\lowerromannumeral}[1]{\romannumeral#1\relax}

\theoremstyle{plain}

\newtheorem{thm}{Theorem}[section]  

\newtheorem{cor}[thm]{{Corollary}} 

\newtheorem{lem}[thm]{{Lemma}}

\newtheorem{prop}[thm]{Proposition}

\newtheorem{defi}[thm]{Definition}

\newtheorem{rmk}[thm]{Remark}

\newcommand{\thistheoremname}{}
\newtheorem*{genericthm*}{\thistheoremname}
\newenvironment{namedthm*}[1]{\renewcommand{\thistheoremname}{#1}%
	\begin{genericthm*}}
	{\end{genericthm*}}

\theoremstyle{remark}

\numberwithin{equation}{section}

\def\Spec{\mathrm{Spec}}

\def\log{\mathrm{log}\,}

\def\Spec{\mathrm{Spec}\,}

\newcommand{\sslash}{\mathbin{/\mkern-6mu/}}

\usepackage[top=1.5in, bottom=1.5in, left=1.5in, right=1.5in]{geometry}

 \usepackage[hyperpageref]{backref}
 
\usepackage{xcolor}
\hypersetup{
	colorlinks,
	linkcolor={red!50!black},
	citecolor={blue!50!black},
	urlcolor={blue!80!black}
}

\contourlength{0.8pt}

\newcommand{\myuline}[1]{%
  \uline{\phantom{#1}}%
  \llap{\contour{white}{#1}}%
}

\begin{document} 
\title[]{Hyperbolicity of varieties with big linear representation of $\pi_1$}

	\author{Ruiran Sun}
		 \address{Department of Mathematics and Statistics, McGill University, Montr\'{e}al, Qu\'{e}bec, Canada, H3A 0B9}
	\email{ruiran.sun@mcgill.ca}

\begin{abstract}
We show the following algebraicity result for a complex projective variety $X$ with big representation of $\pi_1$ into a semi-simple algebraic group: There exists a proper subvariety $Z \subset X$ such that for any algebraic curve $C$, any holomorphic map $\gamma:\, C \to X$ with $\gamma(C) \not\subset Z$ is induced from an algebraic morphism. As an application, we prove pseudo-Brody hyperbolicity of certain varieties with big reductive representations of $\pi_1$.
\end{abstract}


\maketitle

\section{Introduction}
For some complex manifolds it is known that the ``size'' of their fundamental group has influence on geometric properties. In complex dimension one this is clear from the uniformization theorem: a closed Riemann surface is hyperbolic if and only if it has an infinite and non-abelian fundamental group. Things become more complicated in higher dimensions. Hypersurfaces with high degree in projective spaces provide examples of Kobayashi hyperbolic complex manifolds with trivial $\pi_1$. The product of a hyperbolic compact Riemann surface with $\mathbb{P}^n$ is not hyperbolic in any sense, albeit with infinite and non-abelian $\pi_1$.

In his study of the Shafarevich conjecture, Koll\'{a}r introduced the notion of {\it large fundamental group} (cf. \cite{Kol93}), which turns out to be a suitable substitution of infinity of $\pi_1$ in higher dimension. 
A complex projective variety $X$ is said to have large fundamental group if
for every positive dimensional subvariety $Y \subset X$, the image $\mathrm{Im}[\pi_1(Y) \to \pi_1(X)]$ is infinite. People are also interested in varieties with linear $\pi_1$  or linear representations of $\pi_1$, since in this case tools from non-abelian Hodge theory can be applied. In a similar manner, we have
\begin{defi}
Let $G$ be a linear algebraic group defined over $\bar{\mathbb{Q}}$. A Zariski dense representation $\rho:\, \pi_1(X) \to G(\mathbb{C})$ is called {\rm big} if for a sufficiently general point $x \in X$ and every positive dimensional subvariety $Y \subset X$ containing $x$, the image $\rho(\mathrm{Im}[\pi_1(Y) \to \pi_1(X)])$ is infinite.  
\end{defi}
Note that if we take $G$ to be {\it semi-simple}, then the bigness of $\rho$ implies that $\pi_1(X)$ is infinite and non-abelian. Then similar to the complex dimension one case, one can expect {\it hyperbolicity} results for these varieties. In \cite{Zuo96} Zuo proved
\begin{thm}
Suppose that $X$ is a smooth projective variety over $\mathbb{C}$ and $\rho:\,\pi_1(X) \to G(\mathbb{C})$ is a Zariski-dense representation into a semi-simple algebraic group. If $\rho$ is a big representation, then $X$ is Chern hyperbolic, i.e. there exists a proper subvariety $Z \subset X$, such that for any algebraic curve $C$ of genus $g(C) \leq 1$, the image of every non-constant morphism $f:\,C \to X$ is contained in $Z$. 
\end{thm}
Zuo stated the above theorem for almost simple algebraic group $G$ (cf. {\cite[Theorem~2]{Zuo96}}). His argument can be easily generalized to the semi-simple case since one can replace $X$ by a finite \'{e}tale cover and thus $G$ can be assumed to be a direct product of almost simple algebraic groups.

Yamanoi revisited these varieties with Zariski-dense linear representation in \cite{Yam10}, with interest in the {\it distribution of entire curves} (i.e. non-constant holomorphic maps from $\mathbb{C}$ to these varieties).
Motivated by Campana's abelianity conjecture (cf. {\cite[Conjecture~9.8]{Cam04}}),
in {\cite[Proposition~2.1]{Yam10}} Yamanoi proved that for $X$ admitting a Zariski-dense representation of $\pi_1(X)$ into an almost simple algebraic group, every entire curve $\gamma:\,\mathbb{C} \to X$ is degenerate (i.e. the image curve is not Zariski-dense in $X$). A key ingredient of Yamanoi's proof is the combination of the value distribution theory with some construction from non-abelian Hodge theory (i.e. harmonic maps into Bruhat-Tits buildings and spectral coverings,  cf. \S\ref{sec_4.1} and {\cite[\S3]{Yam10}}).

Inspired by Yamanoi's work, we obtain the following algebraicity theorem in this paper: 
\begin{namedthm*}{Main Theorem}\label{main-thm}
Suppose that $X$ is a smooth projective variety over $\mathbb{C}$ and $\rho:\,\pi_1(X) \to G(\mathbb{C})$ is a Zariski-dense representation into a semi-simple algebraic group. If $\rho$ is big, then $X$ is {\rm pseudo Borel hyperbolic}. That is, there exists a proper subvariety $Z \subsetneqq X$ such that for any algebraic curve $C$, any holomorphic map $\gamma:\, C \to X$ with $\gamma(C) \not\subset Z$ is induced from an algebraic morphism. 
\end{namedthm*} 
\begin{rmk}
Note that if $X$ admits such a big representation $\rho$, then for any birational model $\hat{X} \to X$, the pull-back of $\rho$ to $\pi_1(\hat{X})$ is again a big representation. That means the above result cannot be strengthened to $Z = \varnothing$.
\end{rmk}
The reader is referred to \cite{JK20} for the generalities of the notion of Borel hyperbolicity. The above algebraicity result implies pseudo Brody hyperbolicity:
\begin{cor}\label{Cor-1.4}
Let $X$ and $\rho:\,\pi_1(X) \to G(\mathbb{C})$ be the same as in {\rm {\bf Main Theorem}}. If $\rho$ is big, then $X$ is {pseudo Brody hyperbolic}. That is, there exists a proper subvariety $Z \subsetneqq X$ such that the image of every entire curve $\gamma:\,\mathbb{C} \to X$ is contained in $Z$.
\end{cor}
\begin{proof}
Suppose that there exists an entire curve $\gamma:\,\mathbb{C} \to X$ with $\gamma(\mathbb{C}) \not\subset Z$. Note that one can always replace $\gamma$ by the following transcendental holomorphic map
\[
\gamma_{\mathrm{new}}:\, \mathbb{C} \xrightarrow{\mathrm{exp}} \mathbb{C} \xrightarrow{\gamma} X
\]
and we still have $\gamma_{\mathrm{new}}(\mathbb{C}) \not\subset Z$. 

On the other hand, by the Main Theorem we know that $\gamma_{\mathrm{new}}$ is induced from a non-constant morphism $\mathbb{P}^1 \to X$. This gives a contradiction.
\end{proof}

We are now interested in studying Lang's conjecture on varieties with linear fundamental groups. Recall that for a projective variety $X$ the {\it special set} of $X$, $\mathrm{Sp}(X)$, is defined as the Zariski closure of the union of images of all non-constant rational maps from abelian varieties to $X$ (cf. \cite{Lang} or {\cite[\S2.4]{Yam_book}}).
Then Lang's conjecture predicts that $X$ is of general type if and only if $X \neq \mathrm{Sp}(X)$.

In {\cite[Theorem~2.17]{Yam_book}}, Yamanoi considered a normal projective variety $X = S \slash \Gamma$ where $S$ is a Stein space and $\Gamma \subset \mathrm{GL}_n(\mathbb{C})$ is a {\it discrete} linear group which acts freely on $S$.  Motivated by the hyperbolic version of Lang's conjecture, he proved that $X$ is Kobayashi hyperbolic if and only if $\mathrm{Sp}(X) = \varnothing$.

Inspired by Lang's conjecture and Yamanoi's theorem, we show the following
\begin{thm}\label{Yam-question}
Let $X$ be a complex projective variety with a big reductive representation $\rho:\, \pi_1(X) \to \mathrm{GL}_n(\mathbb{C})$. If the special set of $X$ is a proper subset, then $X$ is pseudo-Brody hyperbolic. 
\end{thm}
The strategy of the proof is to consider the abelian part and the semi-simple part of this reductive representation respectively. We use Yamanoi's theorem \cite{Yam15} to deal with the abelian part, and the semi-simple part follows from the Main Theorem of this paper.
Details will be given in Section~\ref{proof-Yam-question}.\\

\noindent{\bf Acknowledgment.} The author is very grateful to Katsutoshi Yamanoi for his generosity in sharing ideas: Theorem~\ref{Yam-question} is previously a question asked by him and the proof given in Section~\ref{proof-Yam-question} relies on his important observation.
The author would also like to thank Steven Lu and Kang Zuo for many useful and inspiring discussions. 

\section{Dichotomy of representations and the strategy of the proof of Main Theorem}

We first investigate the moduli space of representations. Recall that $\mathrm{Hom}(\pi_1(X),G)$ is an affine variety. Then the 
{\it Betti moduli space} is defined as the categorical quotient
\[
M_{\mathrm{B}}(X,G) := \mathrm{Hom}(\pi_1(X),G) \sslash G
\]
which is a quasi-projective scheme over $\bar{\mathbb{Q}}$ (cf. \cite{Sim92}).

A representation $\rho:\, \pi_1(X) \to G(\mathbb{C})$ is said to be {\it rigid} if $[\rho] \in M_{\mathrm{B}}(X,G)(\mathbb{C})$ is an isolated point. Then the rigidity of $\rho$ gives the factorization $\pi_1(X) \to G(K) \subset G(\mathbb{C})$ after conjugation, where $K$ is some number field (cf. \cite{Rag72},p.90, Proposition~6.6).

Let $\mathfrak{p}$ be a prime ideal of the ring of integers of $K$ and $K_{\mathfrak{p}}$ be the $\mathfrak{p}$-adic field. We say $\rho:\,\pi_1(X) \to G(K)$ is {\it $\mathfrak{p}$-bounded} if $\mathrm{Im}[\pi_1(X) \xrightarrow{\rho} G(K) \hookrightarrow G(K_{\mathfrak{p}})]$ is contained in some maximal compact subgroup (see for instance {\cite[\S6]{Zimmer}}).

Following the strategy in \cite{Zuo96}, we consider the dichotomy of representations:
\begin{description}
\item[\myuline{Type A}] \begin{tabular}{l}  
$\rho$ is rigid and the factorization $\pi_1(X) \to G(K)$ is $\mathfrak{p}$-bounded\\ for each $\mathfrak{p} \in \mathrm{Spec}\, \mathcal{O}_K$;\\[.2cm]
\end{tabular}
\item[\myuline{Type B}] 
  \begin{tabular}{l}
    the rest cases, that is, either $\rho:\,\pi_1(X) \to G(\mathbb{C})$ is non-rigid,\\
or $\rho:\, \pi_1(X) \to G(K)$ is $\mathfrak{p}$-unbounded for some $\mathfrak{p} \in \mathrm{Spec}\, \mathcal{O}_K$.
  \end{tabular}
\end{description}

We will see in Section~\ref{Type_A} that the big representation $\rho$ of Type A induces a {variation of Hodge structures} over $X$ with {\it generically injective} Higgs map. In Section~\ref{Type_B} a pluriharmonic map from $\tilde{X}$, the universal covering of $X$, into some Bruhat-Tits building will be constructed from a Type B representation $\rho$. One can consider the {spectral covering} $X^s$ of $X$ associated to this pluriharmonic map, and the bigness of $\rho$ guarantees that $X^s$ has {\it maximal Albanese dimension}.

We will prove the Main Theorem respectively in these two cases.

\section{Type A representations and variations of Hodge structures}\label{Type_A}
In this section we prove the Main Theorem for the case that $\rho$ is a Type A representation.

\subsection{Generalities about Type A representations}

We first recall a lemma about algebraic groups (\cite{Zimmer}, p.120-121).

\begin{lem}
Suppose $H \subset G(K)$ is a subgroup which is $\mathfrak{p}$-bounded for every $\mathfrak{p} \in \Spec \mathcal{O}_K$. Then we have
\[
[H : H \cap G(\mathcal{O}_K)] < +\infty.
\]
\end{lem}

In our case, this means that $\rho(\pi_1(X)) \cap G(\mathcal{O}_K)$ is a finite index subgroup for a Type A representation $\rho$. Therefore, after replacing $X$ by some finite \'{e}tale covering, we can assume further that
\[
\rho:\, \pi_1(X) \to G(\mathcal{O}_K).
\]

Next we shall construct a {\it real semi-simple discrete} representation from $\rho$. We consider the restriction of scalars. Let $\sigma_i:\, K \hookrightarrow \mathbb{C}$, $i=1,2,\dots,d$ be different embeddings. Define
\[
\mathrm{R}_{K/\mathbb{Q}}(G) := \prod^d_{i=1} G^{\sigma_i},
\]
which is an algebraic $\mathbb{Q}$-group with the diagonal embedding
\[
\alpha:\, G(K) \to \mathrm{R}_{K/\mathbb{Q}}(G)(\mathbb{Q}).
\]

Now we consider the composition map
\[
\pi_1(X) \xrightarrow{\rho} G(\mathcal{O}_K) \xrightarrow{\alpha} \mathrm{R}_{K/\mathbb{Q}}(G)(\mathbb{Z}). 
\]

Following Zuo's argument one can find a noncompact factor $G_0$ of the Zariski closure of the image of $\alpha \circ \rho$ in  $\mathrm{R}_{K/\mathbb{Q}}(G)$, such that the induced map
\[
\rho_0:\, \pi_1(X) \to G_0(\mathbb{Z})
\]
is a discrete big representation of $\pi_1(X)$ into $G_0(\mathbb{R})$, a semi-simple real Lie group of noncompact type.

From Simpson correspondence \cite{Sim92}, we know that the rigid representation $\rho$, as well as the induced representation $\rho_0$, have more structures.
First note that the original representation
\[
\rho:\, \pi_1(X) \to G \subset \mathrm{GL}_n(\mathbb{C})
\]
gives us a Higgs bundle $(E,\theta)$ together with a harmonic metric $u$. The rigidity of $\rho$ implies that $(E,\theta)$ is the fixed point of the $\mathbb{C}^*$-action on the moduli space of Higgs bundles. Then by Simpson's ubiquity theorem, we know that $(E,\theta)$ is a Hodge bundle associated to a $\mathbb{C}$-VHS, i.e. there exists a bigrading $E = \bigoplus_{p+q =k} E^{p,q}$ with
\[
\theta|_{E^{p,q}}:\, E^{p,q} \to E^{p-1,q+1} \otimes \Omega^1_X.
\]

For each embedding $\sigma_i:\, K \hookrightarrow \mathbb{C}$, we know that 
\[
\rho^{\sigma_i}:\, \pi_1(X) \to G(K) \xrightarrow{\sigma_i} G(\mathbb{C})
\]
is still rigid, and thus carries a structure of $\mathbb{C}$-VHS by Simpson's theorem. This means that the induced representation
\[
\rho_0:\, \pi_1(X) \to G_0(\mathbb{Z}) \subset \mathrm{R}_{K/\mathbb{Q}}(G)(\mathbb{Z})
\]
is equipped with a structure of $\mathbb{Z}$-VHS.

Denote by $(E_0,\theta_0,u_0)$ the Hodge bundle associated to $\rho_0$ with the harmonic metric. 
By using the theory of harmonic maps and Mok's factorization theorem {\cite[Main Theorem]{Mok92}}, Zuo proved the following
\begin{thm}[{\cite[\S3]{Zuo96}}]\label{VHS}
Suppose $\rho$ is a big representation of Type A. Then the induced representation $\rho_0$ is also big and the Higgs map
\[
\theta_0:\, T_X \to \mathrm{End}(E_0)
 \]
is generically injective.  
\end{thm}

We shall use this construction to prove the pseudo hyperbolicity of $X$.

\subsection{Griffiths line bundle and the big Picard theorem}

Recall that we have a Hodge bundle
\[
\left( E_0 = \bigoplus_{p+q=k}E^{p,q}_0, \theta_0 = \bigoplus_{p+q=k}\theta^{p,q}_0  \right)
\]
coming from a $\mathbb{Z}$-VHS together with a harmonic metric $u_0$ (the Hodge metric).

We consider the {\it Griffiths line bundle} on $X$
\[
\mathrm{K}(E_0) := \left(\mathrm{det}\,E^{k,0}_0\right)^{\otimes k} \otimes \left(\mathrm{det}\,E^{k-1,1}_0\right)^{\otimes (k-1)} \otimes \dots \otimes \left(\mathrm{det}\,E^{1,k-1}_0\right).
\]
Then the curvature form of $\mathrm{K}(E_0)$ with respect to the induced Hodge metric can be written as
\[
\Theta(\mathrm{K}(E_0)) = k\,\mathrm{Tr}\, \Theta(E^{k,0}_0) + (k-1)\,\mathrm{Tr}\, \Theta(E^{k-1,1}_0) + \dots + \mathrm{Tr}\, \Theta(E^{1,k-1}_0). 
\]

\begin{prop}[Griffiths, {\cite[Proposition(7.15)]{Gri70}}] %
For a tangent vector $\eta \in T_X$, one has
\[
\Theta(\mathrm{K}(E_0)) (\eta \wedge \bar{\eta}) \geq 0
\]
with equality if and only if $\theta_{0,\eta} =0$.  
\end{prop}

In our situation, it is easy to show
\begin{lem}\label{lem_Gbig}
The Griffiths line bundle $\mathrm{K}(E_0)$ is big and nef. 
\end{lem}
\begin{proof}
Note that $\mathrm{K}(E_0)$ is nef since $\Theta(\mathrm{K}(E_0))$ is semi-positive. To prove the bigness, one only needs to check that
\[
\int_X \Theta(\mathrm{K}(E_0))^{\mathrm{dim}\,X} >0.
\]
This follows from the fact that $\theta$ is generically injective and thus $\theta_{0,\eta} \neq 0$ for any tangent vector $\eta$ at a general point of $X$.
\end{proof}

Now we can use the {\it Second Main Theorem} of Brotbek-Brunebarbe (cf. {\cite[Theorem~1.1]{BB20}}) to obtain the following
\begin{thm}
Let $(E_0,\theta_0)$ and $\mathrm{K}(E_0)$ be the same as above.
Let $\gamma$ be a holomorphic map from an algebraic curve $C$ to $X$. 
For any ample line bundle $A$ on $X$, there exists a constant $\epsilon>0$ such that
\begin{align}
  \label{SMT-BB}
T(r,\gamma, \mathrm{K}(E_0)) \leq \epsilon \cdot \left(\log r + \log T(r,\gamma,A)  \right),\,\,||.  
\end{align}
\end{thm}
\begin{rmk}
Readers are referred to {\cite[\S2.4]{BB20}} for notations in value distribution theory. Note that in the general form of Brotbek-Brunebarbe's Second Main Theorem, the source space of $\gamma$ is assumed to be a parabolic Riemann surface and the weighted Euler characteristic (cf. {\cite[Definition~1.2]{PS14}}) of it appears in the right hand side of \eqref{SMT-BB}. In our setting, since $C$ is an algebraic curve, we know that the weighted Euler characteristic of $C$ has logarithmic growth (see p.4, (2) of \cite{PS14}).
\end{rmk}

\begin{proof}[Proof of Main Theorem]
From Lemma~\ref{lem_Gbig} we know that the Griffiths line bundle $\mathrm{K}(E_0)$ is big. Then by Kodaira's lemma one can find a positive integer $m$ such that there exists a nonzero section $s$ of $\mathrm{K}(E_0)^{\otimes m} \otimes A^{-1}$.

Now take $Z$ to be the zero locus of $s$. For any holomorphic map $\gamma:\,C \to X$ with $\gamma(C) \not\subset Z$, we apply \eqref{SMT-BB} and obtain
\[
\frac{1}{m} \cdot T(r,\gamma,A) \leq   T(r,\gamma, \mathrm{K}(E_0)) \leq \epsilon \cdot \left(\log r + \log T(r,\gamma,A)  \right),\,\,||.
\]
This means that $T(r,\gamma,A) = O(\log r)$, which implies the algebraicity of $\gamma$ (see {\it e.g.} \cite[2.11.~cas local]{Dem97b} or \cite[Remark 4.7.4.(\lowerromannumeral{2})]{NW}).
\end{proof}

\section{Type B representations and harmonic maps into Bruhat-Tits buildings}\label{Type_B}
In this section we deal with Type B representations.
\subsection{Harmonic maps into buildings and spectral coverings}\label{sec_4.1}
As we mentioned in the introduction, after replacing $X$ by some finite \'{e}tale cover, we can assume that $G \cong G_1 \times \cdots \times G_k$, where $G_i$ are almost simple algebraic groups. We will consider the induced representations $\pi_1(X) \to G \twoheadrightarrow G_i$.

We first consider the $\mathfrak{p}$-unbounded representations. 
By the theory of harmonic maps into Bruhat-Tits buildings due to Gromov and Schoen \cite{GS92}, there exists a non-constant $\rho$-equivariant pluriharmonic map
\[
u_i:\, \tilde{X} \to \Delta(G_i(K_{\mathfrak{p}}))
\]
from the universal covering of $X$ to the Bruhat-Tits building of $G_i(K_{\mathfrak{p}})$ for $i=1,2,\dots,k$.

Next we consider the case that $\rho:\, \pi_1(X) \to G(\mathbb{C})$ is {\it non-rigid}. In this case we know that $M_{\mathrm{B}}(X,G)$ has positive dimensional component and thus one can find an affine curve contained in $M_{\mathrm{B}}(X,G)(\bar{\mathbb{Q}})$. Denote by $T$ the modulo-$p$ reduction of this affine curve to a finite field $k$. 
Choose a compactification $\bar{T}$ and a smooth point $\infty \in \bar{T}\setminus T$. Then we can define the $\infty$-adic valuation $\nu_{\infty}(\bullet)$ on the function field $k(T)$ of $T$, where for any function $f \in k(T)$ the valuation $\nu_{\infty}(f)$ is the vanishing order of $f$ at $\infty$.\\
Now the deformation of representations along $T$ will induce the following representation
\[
\rho_{T,i}:\, \pi_1(X) \to G \left( k(T)_{\infty}\right) \twoheadrightarrow G_i \left( k(T)_{\infty}\right) 
\]
where $k(T)_{\infty}$ is the completion of $k(T)$ under the $\infty$-valuation $\nu_{\infty}(\bullet)$.
Note that $\rho_{T,i}$ is an unbounded representation with respect to the non-archimedean norm induced by the valuation $\nu_{\infty}(\bullet)$. Then again by the theorem of Gromov-Schoen, we can construct a $\rho_{T,i}$-equivariant non-constant pluriharmonic map
\[
u_i:\, \tilde{X} \to \Delta(G_i(k(T)_{\infty}))
\]
from the universal covering to the Bruhat-Tits building of $G_i(k(T)_{\infty})$ for $i=1,2,\dots,k$.\\

Thus in both aforementioned cases of Type B representations, we obtain an equivariant pluriharmonic map $u$ from $\tilde{X}$ to product of Bruhat-Tits buildings ($\prod^k_{i=1}\Delta(G_i(K_{\mathfrak{p}}))$ or $\prod^k_{i=1}\Delta(G_i(k(T)_{\infty}))$). From the complexified differential $\partial u$ one can extract a multi-valued holomorphic one-form $\omega$ on $X$. Following {\cite[\S1, p.146-147]{Zuo96}} we consider a finite ramified Galois covering $\pi:\, X^s \to X$, the {\it spectral covering}, such that $\pi^*\omega$ splits into $l$ single-valued holomorphic one-forms $\omega_1, \dots, \omega_l \in H^0(X^s,\pi^*\Omega^1_X)$. Here $l$ is the number of roots of the Weyl group of the algebraic group $G$. 
Note that the spectral covering $\pi$ is unramified outside the union of zero loci $\bigcup_{i \neq j}(\omega_i - \omega_j)_0$.\\[.1cm]
Now we consider the Albanese map $a:\, X^s \to \mathrm{Alb}(X^s)$. Note that all $\omega_i$'s are pull-back from the Albanese variety. Thus one can find holomorphic one-forms $\tilde{\omega}_i$ on $\mathrm{Alb}(X^s)$ such that $\omega_i = a^*\tilde{\omega}_i$ for $i =1,\dots,l$. Let $B \subset \mathrm{Alb}(X^s)$ be the maximal abelian subvariety such that all $\tilde{\omega}_i$ vanish on it. We set $A:= \mathrm{Alb}(X^s)/B$ and consider the induced morphism
\[
\Phi:\, X^s \xrightarrow{a} \mathrm{Alb}(X^s) \twoheadrightarrow A.
\]
\begin{prop}
If $\rho:\, \pi_1(X) \to G$ is a {\rm big} representation of Type B, then
\begin{itemize}
\item[1)] $\Phi:\, X^s \to A$ is generically finite onto its image.
\item[2)] $X^s$ is of general type. 
\end{itemize}
\end{prop}
\begin{proof}
See \S1 of \cite{Zuo96}.  
\end{proof}
\subsection{Yamanoi's Second Main Theorem and algebraicity}
Now we prove the Main Theorem in the case that $\rho$ is a Type B representation.
First notice that the general case can be easily reduced to the following situation: $C$ is a smooth quasi-projective curve which is a finite ramified covering of $\mathbb{C}$. Denote by $p_C:\, C \to \mathbb{C}$ the covering map. We suppose that $\gamma:\, C \to X$ is a holomorphic map such that the image curve is not contained in the branched locus of $\pi:\, X^s \to X$.

Now we take $Y$ to be the normalization of the fiber product $C \times_X X^s$. Note that a priori $Y$ is only a Riemann surface. Denote by $\tilde{\gamma}:\, Y \to X^s$ the induced holomorphic map. Then we have the following diagram
\[
\xymatrix{
Y \ar[r]^{\tilde{\gamma}} \ar[d]^{\pi_Y} \ar@/_1pc/[dd]_{p_Y} & X^s \ar[r]^{\Phi} \ar[d]^{\pi} & A \\
C \ar[r]^{\gamma} \ar[d]^{p_C} & X & \\
\mathbb{C} & &
}
\]
where the composed map $p_Y:\,Y \to \mathbb{C}$ is a finite surjective holomorphic map (thus $Y$ is a parabolic Riemann surface).

Next we employ tools from Nevanlinna theory (cf. \S3 of \cite{Yam10} for a brief introduction to the notations). For $r>0$, we set $Y(r):= p^{-1}_Y(\mathbb{D}_r)$. Let $\mathrm{Ram}(p_Y)$ be the (analytic) ramification divisor of $p_Y$. We define
\[
N_{\mathrm{ram}\,p_Y}(r) := \frac{1}{\mathrm{deg}\,p_Y} \int^r_0 \#\left( \mathrm{Ram}(p_Y) \cap Y(t)\right) \frac{dt}{t}.
\]
Denote by $R$ the ramification divisor of the spectral covering $\pi$. 
Denote by $R_C$ the ramification divisor of $p_C:\, C \to \mathbb{C}$.
Then according to the pull-back divisors $\tilde{\gamma}^*R$ and $\pi^*_Y R_C$ we have the following partition of $\mathrm{Ram}(p_Y)$:
\[
\mathrm{Ram}(p_Y) = R_1 + R_2
\]
where $R_1 = \tilde{\gamma}^*R$ and $R_2 = \mathrm{Ram}(p_Y) - R_1 $. Since $\pi_Y$ is \'{e}tale outside $R_1$, we know that the ramification over $R_2$ comes from the ramification divisor of $p_C$ and therefore  $R_2 \leq \pi^*_Y R_C$.
 
Thus we have
\[
N_{\mathrm{ram}\,p_Y}(r) = N(r,R_1) + N(r,R_2)
\]
where $N(r,R_i) := \frac{1}{\mathrm{deg}(p_Y)} \int^r_0 \#\left( R_i \cap Y(t)\right) \frac{dt}{t}$. Since $p_C:\, C \to \mathbb{C}$ is algebraic, we have
\[
N(r,R_2) \leq N(r,\pi^*_Y R_C) 
:= \frac{1}{\mathrm{deg}(p_Y)} \int^r_0 \#\left( \pi^*_YR_C \cap Y(t)\right) \frac{dt}{t} =  \frac{\mathrm{deg}(\pi_Y)}{\mathrm{deg}(p_Y)} \int^r_0 \#\left( R_C \cap C(t)\right) \frac{dt}{t} = O(\log r).
\]
Next we want to determine the growth rate of $N(r,\tilde{\gamma},R)$. Let $L$ be an ample line bundle on $X^s$. Yamanoi proved the following
\begin{lem}[]
For any $\varepsilon >0$, we have
\[
N(r,\tilde{\gamma},R) \leq \varepsilon \cdot T(r,\tilde{\gamma},L),\quad ||.
\]
\end{lem}
\begin{proof}
It is known from \cite{Zuo96} that $\mathrm{Supp}\,R \subset \bigcup_{i \neq j}(\omega_i - \omega_j)_0$. In {\cite[p.557, CLAIM]{Yam10}}, Yamanoi proved that for $i \neq j$, 
\[
N(r,\tilde{\gamma},(\omega_i - \omega_j)_0) \leq \varepsilon \cdot T(r,\tilde{\gamma},L), \quad ||.
\]
\end{proof}
Thus we have proved 
\begin{align}\label{N_ram}
N_{\mathrm{ram}\,p_Y}(r) \leq \varepsilon \cdot T(r,\tilde{\gamma},L) + O(\log r),\quad ||.
\end{align}
Now we can apply Yamanoi's Second Main Theorem {\cite[Theorem~1]{Yam15}} for varieties with maximal albanese dimension and obtain the following
\begin{thm}\label{SMT}
Let $L$ be an ample line bundle on $X^s$ and let $\varepsilon$ be a positive constant. Then there exist a proper Zariski closed subset $\Sigma \subsetneqq X^s$ and a positive constant $\alpha$ satisfying the following property: for any holomorphic map $\tilde{\gamma}:\,Y \to X^s$ from any parabolic Riemann surface $p_Y:\, Y \to \mathbb{C}$ such that the image of $\tilde{\gamma}$ is not contained in $\Sigma$, we have
\[
T(r,\tilde{\gamma},K_{X^s}) \leq \alpha \cdot N_{\mathrm{ram}\,p_Y}(r) + \varepsilon \cdot T(r,\tilde{\gamma},L), \quad ||.
\]  
\end{thm}

\begin{proof}[Proof of Main Theorem]
Let $Z$ be the union of $\pi(R)$ and $\pi(\Sigma)$, which is a proper subvariety of $X$. Let $\gamma:\, C \to X$ be a holomorphic map with $\gamma(C) \not\subset Z$. Denote by $\tilde{\gamma}:\, Y \to X^s$ the induced holomorphic map as above.

Recall that $X^s$ is of general type. That means, we can find some positive integer $m$ such that $L \hookrightarrow K^{\otimes m}_{X^s}$. By Theorem~\ref{SMT}, we can find some constant $C>0$ such that
\[
T(r,\tilde{\gamma},L) \leq C \cdot N_{\mathrm{ram}\,p_Y}(r), \quad ||.
\]
Combing this with the previous estimate \eqref{N_ram} of $N_{\mathrm{ram}\,p_Y}(r)$, we finally proved 
\[
T(r,\tilde{\gamma},L) = O(\log r).
\]

Note that we can choose the line bundle $L$ on $X^s$ to be sufficiently positive such that there is a nonzero map $\pi^*H \to L$, where $H$ is an ample line bundle on $X$. Then we have
\begin{align}\label{ineq_H}
T(r,\gamma \circ \pi_Y, H) = T(r, \tilde{\gamma}, \pi^*H) \leq T(r,\tilde{\gamma},L) = O(\log r),
\end{align}
where the first equality comes from the commutative diagram
\[
\xymatrix{
Y \ar[d]^{\pi_Y} \ar[r]^{\tilde{\gamma}} & X^s \ar[d]^{\pi} \\
C \ar[r]^{\gamma} & X
}
\]
By definition of the Nevanlinna order function, we have
\begin{displaymath}
\begin{array}{ll}
T(r,\gamma \circ \pi_Y, H) & = \int^r_0\frac{dt}{t} \int_{\pi^{-1}_Y p^{-1}_C\mathbb{D}(t)} \pi^*_Y\gamma^*c_1(H) \\
& = \int^r_0\frac{dt}{t}\, \mathrm{deg}(\pi_Y) \int_{p^{-1}_C \mathbb{D}(t)} \gamma^*c_1(H) \\
 & = \mathrm{deg}(\pi_Y) \cdot T(r,\gamma,H).
\end{array}
\end{displaymath}
Combining this with \eqref{ineq_H}, we know that $T(r,\gamma,H) = O(\log r)$. Therefore $\gamma:\,C \to X$ is algebraic.
\end{proof}

\section{Proof of Theorem~\ref{Yam-question}}\label{proof-Yam-question}

Denote by $H$ the Zariski closure of the image of the reductive representation $\rho:\,\pi_1(X) \to \mathrm{GL}_n(\mathbb{C})$. After replacing $X$ by some finite etale covering of it, we can assume that $H \cong T \times G_1 \times \cdots \times G_k$, where $T$ is an algebraic torus and $G_i$ are almost simple groups. 
We first show Theorem~\ref{Yam-question} in two extreme cases: the abelian case $H = T$ and the semi-simple case $H = G:= G_1 \times \cdots \times G_k$.

\subsection*{Abelian case $H=T$} We can consider the Stein factorization $X \to Y \to A$ of the Albanese map $X \to A$ induced by the {\it abelian representation} $\rho:\,\pi_1(X) \to T$. Then the bigness of $\rho$ implies that $X \to Y$ is birational. So the assumption on the special set of $X$ also holds on the special set of $Y$. Then we know that $Y$ is a finite covering of an abelian variety with $\mathrm{Sp}(Y) \subsetneqq Y$.
By Kawamata-Ueno fibration theorem, we know that $Y$ is of general type. Then by {\cite[Corollary 1, (1)]{Yam15}} we know that $Y$ (and thus $X$) is pseudo-Brody hyperbolic.

\subsection*{Semi-simple case $H=G$} This case follows from Corollary~\ref{Cor-1.4}.\\

Now we start to prove Theorem~\ref{Yam-question} for the general case (i.e. both the abelian part $T$ and the semi-simple part $G$ are non-trivial). We will consider the induced representations
\[
\rho_T: \pi_1(X) \to H \cong T \times G \twoheadrightarrow T
\]
and 
\[
\rho_G: \pi_1(X) \to H \cong T \times G \twoheadrightarrow G.
\]
By the result of Koll\'{a}r (cf. {\cite[\S3]{Kol93}} or {\cite[\S5.1]{Zuo_book}}), we know that for the representation $\rho$ we have the Shafarevich map
\[
sh_{\rho}:\, X \dashrightarrow \mathrm{Sh}_{\rho}(X)
\]
which is a rational map with connected fibers from $X$ to a normal algebraic variety $\mathrm{Sh}_{\rho}(X)$ (see also \cite{Cam94} for the K\"{a}hler case). Note that the bigness of $\rho$ implies that $sh_{\rho}$ is {\it birational}.

For $\rho_T$ and $\rho_G$, we have the corresponding Shafarevich maps:
\[
sh_{\rho_T}:\, X \dashrightarrow \mathrm{Sh}_{\rho_T}(X), \quad sh_{\rho_G}:\, X \dashrightarrow \mathrm{Sh}_{\rho_G}(X). 
\]
\begin{lem}
The product Shafarevich map
\[
g:=(sh_{\rho_T},sh_{\rho_G}):\, X \dashrightarrow \mathrm{Sh}_{\rho_T}(X) \times \mathrm{Sh}_{\rho_G}(X)
\]
is birational onto its image in $\mathrm{Sh}_{\rho_T}(X) \times \mathrm{Sh}_{\rho_G}(X)$.  
\end{lem}
\begin{proof}
Denote by $W$ the Zariski closure of $g(X)$ in $\mathrm{Sh}_{\rho_T}(X) \times \mathrm{Sh}_{\rho_G}(X)$. Note that Shafarevich maps have connected fibers.
Thus we only need to show that the general fiber of $g:\, X \dashrightarrow W$ is of zero dimension.

Let $F:= g^{-1}(w)$ be a general fiber of $g:\, X \dashrightarrow W$. Since $sh_{\rho_T}(F) = \mathrm{point}$ (resp. $sh_{\rho_G}(F) = \mathrm{point}$), by the property of Shafarevich maps we know that $\mathrm{Im}[\pi_1(F) \to \pi_1(X) \xrightarrow{\rho_T} T]$ (resp. $\mathrm{Im}[\pi_1(F) \to \pi_1(X) \xrightarrow{\rho_G} G]$) is finite. That means $\mathrm{Im}[\pi_1(F) \to \pi_1(X) \xrightarrow{\rho} H]$ is finite. Since $\rho:\, \pi_1(X) \to H \subset \mathrm{GL}_n(\mathbb{C})$ is a big representation, we know that $F$ is of zero dimension. 
\end{proof}
Moreover, since $\rho_T$ is an abelian representation, we know that $\mathrm{Sh}_{\rho_T}(X)$ is an {\it abelian variety} (in fact a quotient of the Albanese variety of $X$).

Now we apply the {\it Factorization theorem} to the semi-simple representation $\rho_G$
(cf. \cite{Mok92} for Type A representations and {\cite[\S4]{Zuo_book}} for Type B representations). Then after replacing $X$ by some finite \'{e}tale covering and blowing up it if necessary, there exists a projective variety $Y$ and a surjective morphism $f:\, X \to Y$ such that $\rho_G:\, \pi_1(X) \to G$ factors through a {\it big} representation $\rho_Y:\,\pi_1(Y) \to G$. By the semi-simple case ($H=G$), we know that $Y$ is pseudo-Brody hyperbolic.

Now we have the following commutative diagram
\[
\xymatrix{
X \ar[d]^f \ar@{-->}[r]^-{g} & \mathrm{Sh}_{\rho_T}(X) \times \mathrm{Sh}_{\rho_G}(X) \ar[d] \\
Y \ar@{-->}[r] & \mathrm{Sh}_{\rho_G}(X)
}
\]
where the general fiber of $Y \dashrightarrow \mathrm{Sh}_{\rho_G}(X)$ is of zero dimension.

Let $Y_0$ be the Zariski open subset of $Y$ where the rational map $Y \dashrightarrow \mathrm{Sh}_{\rho_G}(X)$ is a well-defined morphism.  
Let $X_0:= f^{-1}(Y_0)$. Denote by $\mathscr{A}:= Y_0 \times_{\mathrm{Sh}_{\rho_G}(X)} (\mathrm{Sh}_{\rho_T}(X) \times \mathrm{Sh}_{\rho_G}(X))$ the fiber product. Then $g:\,X_0 \dashrightarrow \mathrm{Sh}_{\rho_T}(X) \times \mathrm{Sh}_{\rho_G}(X)$ factors through $X_0 \dashrightarrow \mathscr{A} \to \mathrm{Sh}_{\rho_T}(X) \times \mathrm{Sh}_{\rho_G}(X)$. Note that $h:\, X_0 \dashrightarrow \mathscr{A}$ is generically finite onto its image in $\mathscr{A}$ since $\mathscr{A} \to \mathrm{Sh}_{\rho_T}(X) \times \mathrm{Sh}_{\rho_G}(X)$ is generically finite and $g$ is birational onto its image.
We have the commutative diagram
\[
\xymatrix{
X_0 \ar[d]^f \ar@{-->}[r]^-{h} & \mathscr{A} \ar[d] \\
Y_0 \ar@{=}[r] & Y_0 
}
\]
such that for a general $y \in Y_0$, the restriction of the birational map $h:\, X_0 \dashrightarrow \mathscr{A}$ on the fiber $F:= f^{-1}(y)$ coincides with the restriction of the Shafarevich map of $\rho_T$ on $F$:
\[
F \hookrightarrow X \stackrel{\rho_T}{\dashrightarrow} \mathrm{Sh}_{\rho_T}(X). 
\]
This restriction of Shafarevich map is induced from the following {\it abelian} representation
\[
\pi_1(F) \to \pi_1(X) \xrightarrow{\rho} H \cong T \times G \twoheadrightarrow T.
\]
Therefore, the restriction $h|_F:\, F \to \mathrm{Sh}_{\rho_T}(X)$ is a well-defined morphism.
We can assume that $F$ is not contained in $\mathrm{Sp}(X)$ by our assumption on the special set. Then by the abelian case ($H=T$), we know that $F$ is also pseudo-Brody hyperbolic.\\

Now we consider an entire curve $\gamma:\, \mathbb{C} \to X$ and the composed holomorphic map $\gamma_Y:\, \mathbb{C} \to X \xrightarrow{f} Y$. Since $Y$ is pseudo-Brody hyperbolic, we know that outside a proper subvariety of $Y$ all the maps $\gamma_Y$ are constant. Thus we only need to consider these entire curves $\gamma$ which are contained in the general fiber $F$. Thus we have
\[
\mathbb{C} \xrightarrow{\gamma} F \xrightarrow{h|_F} \mathrm{Sh}_{\rho_T}(X)
\]
where $h|_F$ is generically finite onto its image in $\mathrm{Sh}_{\rho_T}(X)$.

Now by {\cite[Corollary 1, (2)]{Yam15}}, we know that every entire curve in $F$ is contained in $E\subset F$, where $E$ is the union of $\mathrm{Sp}(F)$ and the exceptional locus of $h|_F$. Note that the special set of $F$ is contained in the special set of $X$ (which is a proper subset by our assumption), and the exceptional locus of $h|_F$ is contained in the exceptional locus of $h$ (which is a proper subset since $h:\,X_0 \dashrightarrow \mathscr{A}$ is generically finite onto its image). Therefore, there exists a proper subset of $X$ such that every entire curve is contained in it.


\begin{thebibliography}{Yam15b}

\bibitem[BB20]{BB20}
Damian Brotbek and Yohan Brunebarbe.
\newblock Arakelov-nevanlinna inequalities for variations of hodge structures
  and applications, 2020.
\newblock arXiv:2007.12957.

\bibitem[Cam94]{Cam94}
Fr\'{e}d\'{e}ric Campana.
\newblock Remarques sur le rev\^{e}tement universel des vari\'{e}t\'{e}s
  k\"{a}hl\'{e}riennes compactes.
\newblock {\it  Bull. Soc. Math. France}, 122(2):255--284, 1994.

\bibitem[Cam04]{Cam04}
Fr\'{e}d\'{e}ric Campana.
\newblock Orbifolds, special varieties and classification theory.
\newblock {\it  Ann. Inst. Fourier (Grenoble)}, 54(3):499--630, 2004.

\bibitem[Dem97]{Dem97b}
Jean-Pierre Demailly.
\newblock Vari\'{e}t\'{e}s hyperboliques et \'{e}quations diff\'{e}rentielles
  alg\'{e}briques.
\newblock {\it  Gaz. Math.}, (73):3--23, 1997.

\bibitem[Gri70]{Gri70}
Phillip~A. Griffiths.
\newblock Periods of integrals on algebraic manifolds. {III}. {S}ome global
  differential-geometric properties of the period mapping.
\newblock {\it  Inst. Hautes \'{E}tudes Sci. Publ. Math.}, (38):125--180, 1970.

\bibitem[GS92]{GS92}
Mikhail Gromov and Richard Schoen.
\newblock Harmonic maps into singular spaces and {$p$}-adic superrigidity for
  lattices in groups of rank one.
\newblock {\it  Inst. Hautes \'{E}tudes Sci. Publ. Math.}, (76):165--246, 1992.

\bibitem[JK20]{JK20}
Ariyan Javanpeykar and Robert Kucharczyk.
\newblock Algebraicity of analytic maps to a hyperbolic variety.
\newblock {\it  Math. Nachr.}, 293(8):1490--1504, 2020.

\bibitem[Kol93]{Kol93}
J\'{a}nos Koll\'{a}r.
\newblock Shafarevich maps and plurigenera of algebraic varieties.
\newblock {\it  Invent. Math.}, 113(1):177--215, 1993.

\bibitem[Lan91]{Lang}
Serge Lang.
\newblock {\it  Number theory. {III}}, volume~60 of {\it  Encyclopaedia of
  Mathematical Sciences}.
\newblock Springer-Verlag, Berlin, 1991.
\newblock Diophantine geometry.

\bibitem[Mok92]{Mok92}
Ngaiming Mok.
\newblock Factorization of semisimple discrete representations of {K}\"{a}hler
  groups.
\newblock {\it  Invent. Math.}, 110(3):557--614, 1992.

\bibitem[NW14]{NW}
Junjiro Noguchi and J\"{o}rg Winkelmann.
\newblock {\it  Nevanlinna theory in several complex variables and {D}iophantine
  approximation}, volume 350 of {\it  Grundlehren der mathematischen
  Wissenschaften [Fundamental Principles of Mathematical Sciences]}.
\newblock Springer, Tokyo, 2014.

\bibitem[PS14]{PS14}
Mihai Paun and Nessim Sibony.
\newblock Value distribution theory for parabolic riemann surfaces, 2014.
\newblock arXiv:1403.6596.

\bibitem[Rag72]{Rag72}
M.~S. Raghunathan.
\newblock {\it  Discrete subgroups of {L}ie groups}.
\newblock Ergebnisse der Mathematik und ihrer Grenzgebiete, Band 68.
  Springer-Verlag, New York-Heidelberg, 1972.

\bibitem[Sim92]{Sim92}
Carlos~T. Simpson.
\newblock Higgs bundles and local systems.
\newblock {\it  Inst. Hautes \'{E}tudes Sci. Publ. Math.}, (75):5--95, 1992.

\bibitem[Yam10]{Yam10}
Katsutoshi Yamanoi.
\newblock On fundamental groups of algebraic varieties and value distribution
  theory.
\newblock {\it  Ann. Inst. Fourier (Grenoble)}, 60(2):551--563, 2010.

\bibitem[Yam15a]{Yam15}
Katsutoshi Yamanoi.
\newblock Holomorphic curves in algebraic varieties of maximal {A}lbanese
  dimension.
\newblock {\it  Internat. J. Math.}, 26(6):1541006, 45, 2015.

\bibitem[Yam15b]{Yam_book}
Katsutoshi Yamanoi.
\newblock Kobayashi hyperbolicity and higher-dimensional {N}evanlinna theory.
\newblock In {\it  Geometry and analysis on manifolds}, volume 308 of {\it 
  Progr. Math.}, pages 209--273. Birkh\"{a}user/Springer, Cham, 2015.

\bibitem[Zim84]{Zimmer}
Robert~J. Zimmer.
\newblock {\it  Ergodic theory and semisimple groups}, volume~81 of {\it 
  Monographs in Mathematics}.
\newblock Birkh\"{a}user Verlag, Basel, 1984.

\bibitem[Zuo96]{Zuo96}
Kang Zuo.
\newblock Kodaira dimension and {C}hern hyperbolicity of the {S}hafarevich maps
  for representations of {$\pi_1$} of compact {K}\"{a}hler manifolds.
\newblock {\it  J. Reine Angew. Math.}, 472:139--156, 1996.

\bibitem[Zuo99]{Zuo_book}
Kang Zuo.
\newblock {\it  Representations of fundamental groups of algebraic varieties},
  volume 1708 of {\it  Lecture Notes in Mathematics}.
\newblock Springer-Verlag, Berlin, 1999.

\end{thebibliography}

\def\cprime{$'$}

\end{document}